\documentclass{amsart}

\usepackage{color}
\usepackage{soul}
\usepackage[pagewise]{lineno}
\usepackage{hyperref}
\hypersetup{
    colorlinks=true, %set true if you want colored links
    linktoc=all,     %set to all if you want both sections and subsections linked
    linkcolor=blue,  %choose some color if you want links to stand out
}

\renewcommand{\phi}{\varphi}

\newcommand{\norm}[1]{\left\| #1 \right\|}
\newcommand{\abs}[1]{\left\vert #1 \right\vert}

\newcommand{\op}[1]{\operatorname{#1}}

\newcommand{\cstar}{$\mbox{C}^*$}

\newcommand{\boldA}{\mathbf{A}}

\newtheorem{theorem}{Theorem}[section]
\newtheorem{lemma}[theorem]{Lemma}

\newtheorem{corollary}[theorem]{Corollary}
\newtheorem{proposition}[theorem]{Proposition}

\newtheorem{fact}[theorem]{Fact}

\theoremstyle{definition}
\newtheorem{definition}[theorem]{Definition}

\numberwithin{equation}{section}

\begin{document}
\title{Computability of states and non-computability of traces}
\author[C. J. Eagle]{Christopher J. Eagle${}^1$} 
\address[C. J. Eagle]
{Department of Mathematics and Statistics, University of Victoria. PO BOX 1700 STN CSC, Victoria, British Columbia, Canada. V8W 2Y2}%
\email{eaglec@uvic.ca}
\urladdr{https://web.uvic.ca/~eaglec/}
\thanks{${}^1$ Supported by NSERC Discovery Grant RGPIN-2021-02459}
\begin{abstract}
We study computable presentations of the state space and trace space of computably presented \cstar-algebras.  We show that the state space has a computably compact presentation with respect to which the evaluation of states is a computable map.  As an application of having computable states, we obtain a computable version of the GNS construction.

We show that the trace space does not share these computational properties.  Specifically, we give an example of a computable presentation of a simple AF algebra for which the trace space has no computably compact presentation with respect to which evaluation is computable.  As a consequence, we demonstrate that the Elliott invariant of a computably presentable classifiable \cstar-algebra is not always presentable in a way that makes all key features of the invariant computable.
\end{abstract}
\maketitle

\section{Introduction}
Building on work of Fox \cite{fox2022}, there has recently been significant activity in studying the extent to which the primary tools of \cstar-algebra theory are effective (in the sense of computable analysis).  For example, Gelfand duality \cite{Gelfand}, \cite{EvaluativePresentations}, functional calculus \cite{FuncCalc}, K-theory \cite{UHFPaper}, \cite{AFPaper}, \cite{Noncomp}, nuclearity \cite{EffectiveNuclearityPaper}, strongly self-absorbing algebras \cite{SSA}, and embedding problems \cite{QWEP}, \cite{MR3532694}, \cite{MR4797476} have all been studied from this point of view.

States and traces have fundamental importance in the theory of \cstar-algebras, especially their representation theory and classification theory.  In this paper we study states and traces from the computability point of view.  

It was shown in \cite{EffectiveNuclearityPaper} that every computable presentation of a unital \cstar-algebra has computable states (see Section \ref{sec:Prelim} below for definitions).  We show that the state space itself has a particularly strong form of computability in Section \ref{sec:State}, and that this in fact characterizes computable presentations of unital \cstar-algebras.  Combining Theorem \ref{thm:PresentationOfStates} and Theorem \ref{thm:StateImpliesComputable}, we have:

\begin{theorem}
Let $\boldA^\#$ be a presentation of a unital \cstar-algebra $\boldA$.  The following are equivalent:
\begin{enumerate}
\item{The presentation $\boldA^\#$ is computable.}
\item{There is a computably compact presentation of the state space $S(\boldA)$ with respect to which the evaluation map $(\phi, a) \mapsto \phi(a)$ is computable.}
\end{enumerate}
\end{theorem}

The proof is uniform, in that indices for the presentation of $S(\boldA)$ and the evaluation map can be computed from a unital index of $\boldA^\#$, and vice versa.  In Section \ref{sec:Repr} we show that from a computable state one can carry out an effective version of the GNS representation.

The Elliott program for classifying \cstar-algebras is based on $K$-theory and traces.  Specifically, given a unital \cstar-algebra $\boldA$, the associated \emph{Elliott invariant} is
\[\op{Ell}(\boldA) = (K_0(\boldA), K_0(\boldA)_+, [1_\boldA]_{K_0}, K_1(\boldA), T(\boldA), \langle \cdot, \cdot \rangle).\]
Here $T(\boldA)$ is the space of traces on $\boldA$ and $\langle \cdot, \cdot \rangle : T(\boldA) \times K_0(\boldA) \to \mathbb{C}$ is defined by $\langle \tau, x \rangle = \tau(x)$.  The culmination of the classification program is the result that simple separable nuclear unital $\mathcal{Z}$-stable \cstar-algebras satisfying the Universal Coefficient Theorem are classified up to isomorphism by their Elliott invariants (see \cite{Winter} for a survey of the classification program).

It is natural to ask whether $\op{Ell}(\boldA)$ can be computed.  In the context of descriptive set theory, it was shown in \cite{Farah} that (in an appropriate coding of the relevant objects as points in standard Borel spaces) the map $\boldA \mapsto \op{Ell}(\boldA)$ is Borel.  To consider the same question from the point of view of computability theory, we must first establish what we mean by a computable presentation of the Elliott invariant. The following expresses the strong form of computability one might hope for from $\op{Ell}(\boldA)$.

\begin{definition}
Let $\boldA$ be a unital \cstar-algebra and let $\boldA^\#$ be a computable presentation of $\boldA$.  A \emph{fully computable presentation of $\operatorname{Ell}(\boldA)$ over $\boldA^\#$} consists of a computable presentation of $(K_0(\boldA), K_0(\boldA)_+, [1_{\boldA}]_{K_0})$, a computable presentation of $K_1(\boldA)$, and a computably compact presentation of $T(\boldA)$ such that, with respect to these presentations:
\begin{itemize}
\item{There is a computable map taking each projection of $M_n(\boldA^\#)$ to its $K_0$ class.}
\item{There is a computable map taking each $K_0$ class to a pair of projections in matrices over $\boldA$ whose difference in $K_0$ is that class.}
\item{There is a computable map taking each unitary in a matrix algebra over $\boldA$ to its $K_1$ class.}
\item{There is a computable map taking each $K_1$ class to a unitary in a matrix algebra over $\boldA$ representing that class.}
\item{The pairing between traces and $K_0$ classes is computable.}
\end{itemize}
\end{definition}

A strong form of proving computability for the Elliott invariant would thus be to show the existence of a Turing functional taking a computable presentation $\boldA^\#$ of a classifiable algebra $\boldA$ to a fully computable presentation of $\op{Ell}(\boldA)$ over $\boldA^\#$.  The main result of Section \ref{sec:Trace} is that such a strong version of computability of the Elliott invariant does not hold.

In \cite{UHFPaper}, the authors construct a computable functor that takes a computably enumerable presentation $\boldA^\#$ of a unital \cstar-algebra to a computably enumerable presentation of $(K_0(\boldA), K_0(\boldA)_+, [1]_{K_0})$ and a computably enumerable presentation of $K_1(\boldA)$, and the $K_0$-functor gives the unique presentation of $K_0$ that satisfies the first two points in the definition above (see Theorem \ref{fact:FindProjection}).  In \cite{Noncomp} an example was given of a computably presentable \cstar-algebra $\boldA$ such that neither $K_0(\boldA)$ nor $K_1(\boldA)$ has any computable presentation at all, and hence $\op{Ell}(\boldA)$ cannot be fully computably presentable over any computable presentation of that algebra.  However, the algebra constructed there is not simple, so it does not fall into the Elliott classification program.  In Section \ref{sec:Trace} we show that there is a \cstar-algebra that satisfies all the regularity properties of the classification theorem but has a computable presentation with respect to which the Elliott invariant has no fully computable presentation.  First, we show:

\begin{theorem}
There is a computable presentation $\boldA^\#$ of a simple AF algebra $\boldA$ such that $T(\boldA)$ does not possess a computably compact presentation making the evaluation pairing computable.
\end{theorem}

Using an effective form of having real rank $0$, we show that this implies the failure of the strong version of computability of the Elliott invariant.

\begin{theorem}
There is a simple AF (in particular, classifiable) algebra $\boldA$ and a computable presentation $\boldA^\#$ of $\boldA$ such that $\op{Ell}(\boldA)$ has no fully computable presentation over $\boldA^\#$.
\end{theorem}

In the Appendix, written with Joaquin Cisneros, we show that the unit is always a computable point of any computable (or, indeed, c.e.) presentation of a unital \cstar-algebra. \\

\noindent \textbf{Acknowledgements}: The author is grateful to Joaquin Cisneros, Isaac Goldbring, Timothy McNicholl, and Russell Miller for helpful discussions on the topics of this paper. \\
\noindent \textbf{AI Disclosure}: ChatGPT 5.6 was used to assist with literature search, exploring early ideas, and proofreading.  The arguments presented here were written solely by the authors, who take full responsibility for their content.

\section{Preliminaries}\label{sec:Prelim}
We use the framework of computability theory for \cstar-algebras that is presented in \cite{EMST}, and refer the reader there for an introduction to the subject.

To avoid repetition in the statements of theorems, \emph{throughout the paper $\boldA$ is a separable unital \cstar-algebra and $\boldA^\#$ is a fixed presentation of $\boldA$}.  Except where explicitly noted otherwise, we assume that $\boldA^\#$ is a computable presentation.  Whenever we consider $\mathbb{C}^k$ we equip it with the standard presentation whose special points are $\mathbb{Q}(i)^k$. 

We will make use of the fact that $1_\boldA$ is a computable point of $\boldA^\#$; see the Appendix for the proof that this is always the case.  A \emph{unital index} for $\boldA^\#$ is $\langle e_1, e_2\rangle$, where $e_1$ is an index for $\boldA^\#$, $e_2$ is an $\boldA^\#$-index for $1_\boldA$, and $\langle \cdot, \cdot \rangle$ is the standard computable bijection from $\mathbb{N}^2$ to $\mathbb{N}$.

We briefly recall some specific notions from \cstar-algebras and computable analysis that will play a prominent role in the paper.

\begin{definition}
A \emph{state} on $\boldA$ is a linear functional $\omega : \boldA \to \mathbb{C}$ that satisfies $\omega(1_\boldA) = 1$ and $\omega(a^*a) \geq 0$ for all $a \in \boldA$.  A state $\omega$ is \emph{tracial} if it additionally satisfies $\omega(ab) = \omega(ba)$ for all $a, b \in \boldA$.  We denote the space of states on $\boldA$ by $S(\boldA)$ and the space of tracial states by $T(\boldA)$.
\end{definition}

There are \cstar-algebras with $T(\boldA) = \emptyset$ (for example, this occurs in purely infinite \cstar-algebras, such as the Cuntz algebras $\mathcal{O}_n$).  By contrast, every \cstar-algebra has many states (see, for example, \cite[Corollary 3.3.6]{Murphy}).  Every state is bounded and in fact has norm $1$, because for any positive functional $\rho$ we have $\norm{\rho} = \rho(1)$ \cite[Corollary 3.3.4]{Murphy}.  It follows from the Banach-Alaoglu theorem that both $S(\boldA)$ and $T(\boldA)$ are compact and metrizable in the weak${}^*$-topology on the dual of $\boldA$.  When we consider $S(\boldA)$ and $T(\boldA)$ it is always with the weak${}^*$-topology.  It is immediate from the definitions that both $S(\boldA)$ and $T(\boldA)$ are convex.

The following definition was introduced in \cite{Yasgui.1999}.
\begin{definition}
A computable presentation $X^\# = (X, d, (p_n)_{n \in \mathbb{N}})$ of a compact Polish space $(X, d)$ is \emph{computably compact} if there is an algorithm that accepts $k \in \mathbb{N}$ and outputs $n_0, \ldots, n_t \in \mathbb{N}$ such that $X = \bigcup_{j \leq t}B(p_{n_j}; 2^{-k})$.
\end{definition}
Computable compactness is an extremely robust notion and has been shown to be the natural notion of compactness in computable analysis.  For a detailed treatment of computable compactness, see \cite{ComputablyCompact}.  As one example of the naturality of computable compactness, in \cite{Gelfand} it was shown that Gelfand duality holds between computable presentations of commutative unital \cstar-algebras and computably compact presentations of compact metric spaces.

\begin{definition}
Let $X^\# = (X, d, (p_n))_{n \in \mathbb{N}}$ be a computable presentation of a Polish space $X$.  A set $A \subseteq X$ is:
\begin{itemize}
    \item{\emph{c.e. open} if there is a c.e. set $I \subseteq \mathbb{N}^2$ such that $A = \bigcup_{(i,j) \in I}B(p_i; 2^{-j})$.}
    \item{\emph{co-c.e. closed} if $X \setminus A$ is c.e. open}
    \item{\emph{c.e. closed} if $\{(i, j) \in \mathbb{N}^2 : B(p_i; 2^{-j}) \cap A \neq \emptyset\}$ is c.e..}
\end{itemize}
\end{definition}

In the literature, co-c.e. closed sets are sometimes referred to as being \emph{effectively closed} and c.e. closed sets are also known as \emph{effectively overt}.  

A key fact about c.e. closed sets is that if $C \subseteq X$ is c.e. closed in $X^\#$ then $C$ contains a computable dense sequence (see \cite[Lemma 3.27]{ComputablyCompact}), and hence inherits a computable presentation of its own (using the computable dense sequence as special points and the restriction of the metric from $X^\#$).  When the ambient space $X$ is computably compact, the induced presentation on a c.e. closed set $C$ is computably compact if and only if $C$ is also co-c.e. closed (see \cite[Proposition 3.29]{ComputablyCompact}).

When considering $S(\boldA)$ or $T(\boldA)$, it is natural to want not only a computably compact presentation, but additionally one that allows us to compute the evaluation of a state on a point of $\boldA$.  We make the following definition, which is modelled after \cite[Definition 3.1]{EvaluativePresentations}.

\begin{definition}
A presentation $S(\boldA)^\#$ of $S(\boldA)$ is \emph{evaluative} over $\boldA^\#$ if the evaluation map $ev : S(\boldA) \times \boldA \to \mathbb{C}$ given by $\op{ev}(\omega, a) = \omega(a)$ is a computable map from the presentation of $S(\boldA) \times \boldA$ induced by $S(\boldA)^\#$ and $\boldA^\#$.  The definition of an evaluative presentation of $T(\boldA)$ is analogous.
\end{definition}

\section{Computability of the state space}\label{sec:State}
For ease of exposition, we make some definitions that will be in place for the entirety of this section.  We have our fixed computable presentation $\boldA^\#$ of a unital \cstar-algebra $\boldA$.  Let $(a_n)_{n \in \mathbb{N}}$ be an effective enumeration of the generated points of $\boldA^\#$.  For each $n \in \mathbb{N}$, let $\mathbb{D}_n = \{z \in \mathbb{C} : \abs{z} \leq \norm{a_n}\}$.  Let $X = \prod_{n \in \mathbb{N}}\mathbb{D}_n$, viewed with the product topology, and let $\Phi : S(\boldA) \to X$ be defined by
\[\Phi(\omega) = (\omega(a_n))_{n \in \mathbb{N}}.\]
We equip $S(\boldA)$ with the metric
\[d_S(\omega, \phi) = \sum_{n=0}^\infty 2^{-n-1}\min\{1, \abs{\omega(a_n) - \phi(a_n)}\},\]
which (because $(a_n)_{n \in \mathbb{N}}$ is dense in $\boldA$) is compatible with the weak${}^*$-topology on $S(\boldA)$.  Likewise, we equip $X$ with the metric
\[d_X(x, y) = \sum_{n=0}^{\infty}2^{-n-1}\min\{1, \abs{x_n-y_n}\},\]
which is compatible with the product topology.  Note that, with these metrics, $\Phi$ is an isometric embedding of $S(\boldA)$ into $X$.  Let $X^\#$ be the computable presentation of $(X, d_X)$ whose special points are sequences of the form $(x_n)_{n \in \mathbb{N}}$ such that each $x_n \in \mathbb{Q}(i)$ and satisfies either $x_n = 0$ or $\abs{x_n} < \norm{a_n}$, and such that only finitely many coordinates are non-zero.

\begin{lemma}\label{lemma:StatesCoCE}
The image of $\Phi$ is co-c.e. closed in $X^\#$.
\end{lemma}
\begin{proof}
We will describe the image of $\Phi$ as the intersection of three co-c.e. closed sets, corresponding to the three defining properties of states.

First, consider linearity.  Let $(q_k)_{k \in \mathbb{N}}$ be an effective enumeration of $\mathbb{Q}(i)$.  Since the algebraic operations are computable on the generated points of $\boldA^\#$, there is a computable function $f : \mathbb{N}^4 \to \mathbb{N}$ such that for all $(i, j, k, l) \in \mathbb{N}^4$, $a_{f(i,j,k,l)} = q_ka_i + q_la_j$.  For each $(i, j, k, l) \in \mathbb{N}^4$, define
\[L_{i,j,k,l} = \{x \in X : x_{f(i,j,k,l)} = q_kx_i + q_lx_j\}.\]
Since $f$ is computable and the coordinate projections are uniformly computable, $L_{i,j,k,l}$ is co-c.e. closed uniformly in $i,j,k,l$.  Let
\[L = \bigcap_{(i,j,k,l) \in \mathbb{N}^4}L_{i,j,k,l}.\]
Then $L$ is co-c.e. closed.

Next, consider positivity.  In the same way as above, there is a computable function $g : \mathbb{N} \to \mathbb{N}$ such that for each $k \in \mathbb{N}$, we have $a_{g(k)} = a_k^*a_k$.  For each $k \in \mathbb{N}$, define
\[P_k = \{x \in X : x_{g(k)} \in [0, \Vert a_{g(k)} \Vert]\}.\]
Since $g$, $\norm{a_{g(k)}}$, and the coordinate projections are computable uniformly in $k$, the sets $P_k$ are co-c.e. closed uniformly in $k$.  Let
\[P = \bigcap_{k \in \mathbb{N}}P_k.\]
Then $P$ is co-c.e. closed.

Finally, we consider the property of being unital.  Since $1_\boldA$ is a computable point of $\boldA^\#$, there is a computable function $h : \mathbb{N} \to \mathbb{N}$ such that $(a_{h(n)})_{n \in \mathbb{N}}$ is a name for $1_\boldA$.  For each $n \in \mathbb{N}$, let
\[U_n = \{x \in X : \abs{x_{h(n)} - 1} \leq 2^{-n}\}.\]
As in the previous cases, the sets $U_n$ are co-c.e. closed, uniformly in $n$.  Let
\[U = \bigcap_{n \in \mathbb{N}}U_n.\]

We have established that $L \cap P \cap U$ is co-c.e. closed.  The fact that $L \cap P \cap U \supseteq \Phi(S(\boldA))$ follows immediately from the definitions, while the other inclusion follows because any element of $L \cap P$ defines a $\mathbb{Q}(i)$-linear positive map on the dense $\mathbb{Q}(i)$-*-subalgebra of $\boldA$ consisting of generated points, which extends uniquely by continuity to a positive linear functional on $\boldA$, while the condition in $U$ ensures that this functional is unital.
\end{proof}

We next show that the image of $\Phi$ is c.e. closed in $X^\#$.  To do this, we will make use of results about computable states obtained in \cite{EffectiveNuclearityPaper} and transport them to $X^\#$ via the following lemma.

\begin{lemma}\label{lem:MoveComputable}
Let $\omega$ be a computable state on $\boldA^\#$.  Then $\Phi(\omega)$ is a computable point of $X^\#$.  Moreover, from an index for $\boldA^\#$ and an index for $\omega$ we can compute an index for $\Phi(\omega)$.
\end{lemma}
\begin{proof}
Given $k$, compute $N$ large enough so that $\sum_{n > N}2^{-n-1} < 2^{-k-1}$.  Compute a rational $\epsilon$ such that $0 < \epsilon < \frac{2^{-k-1}}{\sum_{n=0}^N2^{-n-1}}$.  

For each $n \in \{0, \ldots, N\}$, we compute a point in $\mathbb{D}_n$ that is within $\epsilon$ of $\omega(a_n)$, as follows.  Given $n$, since $\norm{a_n}$ is computable, we can (in parallel) search for $q_0 \in \mathbb{Q}$ such that $\norm{a_n} < q_0 < \epsilon/2$ and also search for $q_1 \in \mathbb{Q}$ such that $0 < q_1$ and $\norm{a_n} - \epsilon/8 < q_1 < \norm{a_n}$.  At least one of these searches will succeed.  If the former search succeeds first or both searches succeed in the same number of steps, set $r_n = 0$.  In this case, we have $r_n \in \mathbb{D}_n$ and 
\[\abs{r_n - \omega(a_n)} = \abs{\omega(a_n)} \leq \norm{a_n} < q_0 < \epsilon/2 < \epsilon.\]
Now suppose that the second search succeeds first.  Using that $\omega$ is computable, compute $z \in \mathbb{Q}(i)$ such that $\abs{z - \omega(a_n)} < \epsilon/8$.  Compute also $p \in \mathbb{Q}$ such that $\abs{z} < p < \abs{z}+\epsilon/8$.  Define
\[r_n = \frac{q_1}{\max\{q_1, p\}}z.\]
Note that $r_n \in \mathbb{Q}(i)$, and
\[\abs{r_n} = \frac{q_1}{\max\{q_1, p\}}\abs{z} < q_1 < \norm{a_n},\]
so $r_n \in \mathbb{D}_n$.
If $p \leq q_1$ then $r_n = z$, so 
\[\abs{r_n - \omega(a_n)} = \abs{z-\omega(a_n)} < \epsilon/8 < \epsilon.\]
If $p > q_1$ then:
\begin{align*}
\abs{r_n - z} &= \abs{\frac{q_1}{p}z - z} \\
&= \abs{1 - \frac{q_1}{p}}\abs{z} \\
&= \left(\frac{p - q_1}{p}\right)\abs{z} \\
&< p - q_1 \\
&= (p-\abs{z}) + (\abs{z}-\norm{a_n}) + (\norm{a_n}-q_1) \\
&\leq \frac{\epsilon}{8} + (\abs{z-\omega(a_n)} + \abs{\omega(a_n)} - \norm{a_n}) + \frac{\epsilon}{8} \\
&\leq \frac{\epsilon}{8} + (\frac{\epsilon}{8} + \norm{a_n}-\norm{a_n}) + \frac{\epsilon}{8} \\
&< \frac{\epsilon}{2}
\end{align*}
Thus
\[\abs{r_n - \omega(a_n)} \leq \abs{r_n - z} + \abs{z - \omega(a_n)} < \frac{\epsilon}{2} + \frac{\epsilon}{8} < \epsilon.\]

Now define $x_k = (r_0, r_1, \ldots, r_N, 0, 0, \ldots) \in X$.  Then we have
\begin{align*}
d_X(x_k, \Phi(\omega)) &= \sum_{n=0}^{\infty}2^{-n-1}\min\{1, \abs{(x_k)_n - \omega(a_n)}\} \\
&= \sum_{n=0}^N2^{-n-1}\min\{1, \abs{r_n - \omega(a_n)}\} + \sum_{n > N}2^{-n-1}\min\{1, \abs{\omega(a_n)}\} \\
&\leq \epsilon\sum_{n=0}^N2^{-n-1} + \sum_{n>N}2^{-n-1} \\
&< 2^{-k-1}+2^{-k-1} \\
&= 2^{-k}
\end{align*}
Thus $(x_k)_{k\in\mathbb{N}}$ is a computable name for $\Phi(\omega)$.  The second claim follows from the uniformity of the proof.
\end{proof}

We next need a uniform version of \cite[Corollary 2.9(3)]{EffectiveNuclearityPaper}.  This uniform version follows in the same way as the version stated in \cite{EffectiveNuclearityPaper}, but since a proof was not given in \cite{EffectiveNuclearityPaper} we include one here for completeness.

Given a finite tuple $\vec{c} = (c_1, \ldots, c_k) \in \boldA^k$, the \emph{joint numerical range} of the tuple is
\[W(\vec{c}) = \{(\omega(c_1), \ldots, \omega(c_k)) \in \mathbb{C}^k : \omega \in S(\boldA)\}.\]
The joint numerical range of a tuple of generated points can be approximated by computable states, in the following sense:
\begin{fact}[{\cite[Corollary 2.9(2)]{EffectiveNuclearityPaper}}]\label{fact:EffNuc}
There is an algorithm that, given generated points $\vec{c} = (c_1, \ldots, c_k) \in \boldA^k$, a tuple $\vec{b} = (b_1, \ldots, b_k) \in \mathbb{Q}(i)^k$, and a rational $\delta > 0$ such that $d(\vec{b}, W(\vec{c})) < \delta$, returns the index of a computable state on $\boldA^\#$ such that $\max_{1 \leq i \leq k}\abs{b_i - \omega(c_i)} < \delta$.  Here $d$ refers to the $\ell_\infty$ metric on $\mathbb{C}^k$.
\end{fact}

\begin{lemma}\label{lem:densesequencestates}
There is a computable sequence of states on $\boldA^\#$ that is weak${}^*$-dense in $S(\boldA)$.  Moreover, an index for such a sequence can be computed from a unital index for $\boldA^\#$.
\end{lemma}
\begin{proof}
We describe an algorithm for producing such a sequence.  

Let $(\vec{c}_n)_{n \in \mathbb{N}}$ be an effective listing of all finite tuples of generated points of $\boldA^\#$, let $(\vec{b}_n)_{n \in \mathbb{N}}$ be an effective listing of all finite tuples from $\mathbb{Q}(i)$, and $(\delta_n)_{n \in \mathbb{N}}$ be an effective listing of the positive rational numbers.  Apply Fact \ref{fact:EffNuc} to the generated point $1$, the rational number $1$, and the error $\delta = 1$ to obtain an index for a computable state $\phi$.

Given $n$, express $n = \langle m, t, r, s \rangle$, where $\langle \rangle$ is the standard computable bijection between $\mathbb{N}^4$ and $\mathbb{N}$.  If the length of the tuple $\vec{c}_m$ is not the same as the length of the tuple $\vec{b}_t$, set $\omega_n = \phi$.

Now, assuming that $\vec{c}_m$ and $\vec{b}_t$ are of the same length, by \cite[Corollary 2.9(1)]{EffectiveNuclearityPaper} the set $W(\vec{c}_m)$ is a computably compact subset of $\mathbb{C}^{k_m}$ (where $k_m$ is the length of $\vec{c}_m$), and we can compute an index for $W(\vec{c}_m)$.  From that index, using computable compactness, it follows that $d(\vec{b}_t, W(\vec{c}_m))$ is computable (uniformly in $t$ and $m$).  We can thus compute a rational $q_{m,t,s}$ such that $q_{m,t,s} \leq d(\vec{b}_t, W(\vec{c}_m)) < q_{m,t,s}+2^{-s}$.  Now we can test whether or not $q_{m,t,s}+2^{-s} < \delta_r$.  If not, set $\omega_n = \phi$.  If so, set $\omega_n$ to be the state computed from $\vec{c}_m$, $\vec{b}_t$, and $\delta_r$ by Fact \ref{fact:EffNuc}.

Finally, to see that $(\omega_n)_{n \in \mathbb{N}}$ is weak${}^*$-dense, suppose we have a state $\omega \in S(\boldA)$ and a basic weak${}^*$-neighbourhood $U$ of $\omega$.  Then, by density, there are generated points $\vec{c} = (c_1, \ldots, c_k)$ of $\boldA^\#$ and a rational $\epsilon > 0$ such that
\[\{\psi \in S(\boldA) : \max_{1 \leq i \leq k}\abs{\psi(c_i) - \omega(c_i)} < \epsilon\} \subseteq U.\]
Choose $\vec{b} = (b_1, \ldots, b_k) \in \mathbb{Q}(i)^k$ such that 
\[\max_{1 \leq i \leq k}\abs{b_i - \omega(c_i)} < \frac{\epsilon}{4}.\]\
Then
\[d(\vec{b}, W(\vec{c})) \leq d(\vec{b}, (\omega(c_1), \ldots, \omega(c_k))) = \max_{1 \leq i \leq k}\abs{b_i - \omega(c_i)} < \frac{\epsilon}{4}.\]
For some $m, t, r, s$, we have $\vec{c} = \vec{c}_m$, $\vec{b} = \vec{b}_t$, $\epsilon/2 = \delta_r$, and $2^{-s} < \epsilon/4$.  Let $n = \langle m, t, r , s\rangle$.  With these parameters we have \[q_{m,t,s} + 2^{-s} \leq d(\vec{b}_t, W(\vec{c}_m)) + 2^{-s} < \epsilon/4+\epsilon/4 = \epsilon/2 = \delta_r.\]
Thus, by the construction above, we have
\[\max_{1 \leq i \leq k}\abs{b_i - \omega_n(c_i)} < \epsilon/2.\]
And so
\begin{align*}
\max_{1 \leq i \leq k}\abs{\omega_n(c_i) - \omega(c_i)} &\leq \max_{1 \leq i \leq k}\abs{\omega_n(c_i) - b_i} + \max_{1 \leq i \leq k}\abs{\omega(c_i) - b_i} \\
&< \epsilon/2 +  \epsilon/4 \\
&< \epsilon.
\end{align*}
Thus $\omega_n \in U$.
\end{proof}

\begin{theorem}\label{thm:PresentationOfStates}
Let $\boldA^\#$ be a computable presentation of a unital \cstar-algebra $\boldA$.  There is a computably compact evaluative presentation of $S(\boldA)$ over $\boldA^\#$.  Moreover, an index for such a presentation can be computed from a unital index for $\boldA^\#$.
\end{theorem}
\begin{proof}
Using Lemma \ref{lem:densesequencestates}, we can compute an index for a computable sequence $(\omega_n)_{n \in \mathbb{N}}$ of states on $\boldA$ that is weak${}^*$-dense in $S(\boldA)$.  By Lemma \ref{lem:MoveComputable}, the sequence $(\Phi(\omega_n))_{n \in \mathbb{N}}$ is a computable dense sequence in $X^\#$.  Thus, by \cite[Lemma 3.27]{ComputablyCompact}, the image of $\Phi$ is c.e. closed in $X^\#$.  In Lemma \ref{lemma:StatesCoCE} we showed that the image of $\Phi$ is co-c.e. closed in $X^\#$.  Together, by \cite[Proposition 3.29]{ComputablyCompact}, these facts imply that the presentation $(\Phi(S(\boldA)), (\Phi(\omega_n))_{n \in \mathbb{N}}, d_X|_{\Phi(S(\boldA))})$ is a computably compact presentation of the image of $\Phi$.  Since $\Phi$ is an isometry, this implies that $S(\boldA)^\# = (S(\boldA), (\omega_n)_{n \in \mathbb{N}}, d_S)$ is a computably compact presentation of $S(\boldA)$.

It remains to verify that $S(\boldA)^\#$ is evaluative over $\boldA^\#$.  We are given an $\boldA^\#$-name of a point $a \in \boldA$ and an $S(\boldA)^\#$-name of an $\omega \in S(\boldA)$ and a $k \in \mathbb{N}$.  Using the name for $a$, compute $n$ such that $\norm{a-a_n} < 2^{-k-2}$.  Likewise, using the name for $\omega$, compute $m$ such that $d_S(\omega, \omega_m) < 2^{-k-n-4}$.  Finally, compute $q \in \mathbb{Q}(i)$ such that $\abs{\omega_m(a_n)-q} < 2^{-k-2}$.  We will show that $\abs{\op{ev}(\omega, a) - q} < 2^{-k}$.

By definition of $d_S$, we have
\[d_S(\omega, \omega_m) \geq 2^{-n-1}\min\{1, \abs{\omega(a_n) - \omega_m(a_n)}\},\]
and thus,
\begin{align*}
\min\{1, \abs{\omega(a_n) - \omega_m(a_n)}\} &\leq 2^{n+1}d_S(\omega, \omega_m) \\
&< 2^{n+1}2^{-k-n-4} \\
&=2^{-k-3}.
\end{align*}
Since $2^{-k-3} < 1$, this implies
\[\abs{\omega(a_n) - \omega_m(a_n)} < 2^{-k-3}.\]
Thus
\begin{align*}
\abs{\omega(a) - q} 
&< \abs{\omega(a) - \omega(a_n)} + \abs{\omega(a_n) - \omega_m(a_n)} + \abs{\omega_m(a_n) - q} \\
&\leq \abs{\omega(a-a_n)} + 2^{-k-3} + 2^{-k-2} \\
&\leq \norm{a-a_n}+2^{-k-3}+2^{-k-2} \\
&< 2^{-k-2}+2^{-k-3}+2^{-k-2} \\
&< 2^{-k}.
\end{align*}
\end{proof}

We also have the converse:

\begin{theorem}\label{thm:StateImpliesComputable}
Suppose that $\boldA^\dagger$ is a presentation of $\boldA$ and that $S(\boldA)$ has a presentation $S(\boldA)^\dagger$ that is computably compact and evaluative over $\boldA^\dagger$.  Then $\boldA^\dagger$ is computable.  Moreover, a unital index for $\boldA^\dagger$ can be computed from an index for $S(\boldA)^\dagger$ and an index for the evaluation map.
\end{theorem}
\begin{proof}
For any $a \in \boldA$ and any $\omega \in S(\boldA)$ we have $\abs{\omega(a)} \leq \norm{\omega}\norm{a} = \norm{a}$.  If $a$ is normal, then there is some state $\omega_a$ such that $\norm{a} = \abs{\omega(a)}$ \cite[Theorem 3.3.6]{Murphy}.  Thus, if $a$ is normal, we have
\[\norm{a} = \sup_{\omega \in S(\boldA)}\abs{\omega(a)}.\]
Now if $a$ is a generated point of $\boldA^\dagger$, define $f : S(\boldA) \to \mathbb{R}$ by 
\[f(\omega) = \abs{\omega(a^*a)} = \abs{\operatorname{ev}(\omega, a^*a)}.\]
Since $S(\boldA)^\dagger$ is evaluative over $\boldA^\dagger$, $f$ is computable, uniformly in $a$.  We have
\[\norm{a} = \sqrt{\norm{a^*a}} = \sqrt{\sup_{\omega \in S(\boldA)}\abs{\omega(a^*a)}} = \sqrt{\sup_{\omega \in S(\boldA)}f(\omega)}.\]
Since $S(\boldA)^\dagger$ is computably compact and $f$ is computable, the supremum is computable, uniformly in $f$ (see \cite[Proposition 3.5(1)]{ComputablyCompact}).  Thus $\norm{a}$ is computable, uniformly in $a$.

To compute an index for the unit, let $(b_n)_{n \in \mathbb{N}}$ be an effective enumeration of the generated points of $\boldA^\dagger$.  For each $n$, we have
\begin{align*}
\norm{1_\boldA - b_n}^2 &= \norm{(1_\boldA - b_n)^*(1_\boldA-b_n)} \\
&= \sup_{\omega \in S(\boldA)}\omega((1_\boldA - b_n)^*(1_\boldA-b_n)) \\
&= \sup_{\omega \in S(\boldA)}(\omega(1_\boldA) - \omega(b_n^*) - \omega(b_n) + \omega(b_n^*b_n)) \\
&= \sup_{\omega \in S(\boldA)}(1 - \omega(b_n^*) - \omega(b_n) + \omega(b_n^*b_n))
\end{align*}
Since the algebraic operations are computable on the generated points, given $n$ we can compute $r$ and $s$ such that $b_n^* = b_r$ and $b_n^*b_n = b_s$.  Thus, by evaluativity, $1 - \omega(b_n^*)-\omega(b_n)+\omega(b_n^*b_n)$ is computable, uniformly in $n$ and $\omega$.  From the calculation above, it follows as in the proof that $\norm{a}$ is computable that $\norm{1_\boldA - b_n}$ is computable, uniformly in $n$.  To build a computable name for $1_\boldA$, given $k$, search to find an $n(k)$ such that $\norm{1_\boldA - b_{n(k)}} < 2^{-k}$.  Since the generated points are dense, this search terminates.  Then $(b_{n(k)})_{k \in \mathbb{N}}$ is a computable name for $1_\boldA$.
\end{proof}

In general, we do not know whether every computable presentation of a unital \cstar-algebra admits a computable pure state.  For commutative algebras this is so, because pure states correspond to point evaluations on the Gelfand spectrum, so the existence of computable pure states follows from \cite[Theorem 1.1]{EvaluativePresentations}.

\section{Effective GNS representations}\label{sec:Repr}
A primary use of states is to create representations of \cstar-algebras on Hilbert spaces, through the GNS construction.  We show that this can be carried out effectively.  In this section we only require $\boldA^\#$ to be c.e., so we explicitly state this hypothesis in each theorem.

We begin by setting some notation and recalling the classical GNS construction (for more details, see, for example, \cite[Section 3.4]{Murphy}).  Suppose that $\phi$ is a state on $\boldA$.  Let $N_\phi = \{a \in \boldA : \phi(a^*a) = 0\}$.  Let $\langle \cdot, \cdot \rangle_\phi : \boldA \times \boldA \to \mathbb{C}$ be defined by $\langle a, b \rangle_\phi = \phi(b^*a)$.  This is a semi-inner product, which induces a semi-norm on $\boldA$ by $\norm{a}_\phi = \sqrt{\langle a, a \rangle_\phi}$.  The completion of $\boldA/N_\phi$ is thus a Hilbert space $H_\phi$, whose inner product we also denote by $\langle \cdot, \cdot \rangle_\phi$ and whose norm we denote by $\norm{\cdot}_\phi$.  We denote by $\eta_\phi : \boldA \to H_\phi$ the canonical map given by sending $a$ to $a+N_\phi$ in $\boldA/N_\phi$ and then following the standard embedding of $\boldA/N_\phi$ into $H_\phi$.  Note that $\norm{\eta_\phi} \leq 1$.

\begin{definition}
Let $\phi$ be a state on $\boldA$.  A \emph{presentation of $H_\phi$ induced by $\boldA^\#$} is a presentation that has special points $(\eta_\phi(a_j))_{j \in \mathbb{N}}$, where $(a_j)_{j \in \mathbb{N}}$ is an effective enumeration of the generated points of $\boldA^\#$.
\end{definition}

We note that since $(a_j)_{j \in \mathbb{N}}$ is dense in $A$ it follows that $(\eta_\phi(a_j))_{j \in \mathbb{N}}$ is dense in $H_\phi$, so a presentation of $H_\phi$ induced by $\boldA^\#$ really is a presentation.  Determining membership in $N_\phi$ may not be effectively possible, so we do not attempt to make the list $(\eta_\phi(a_j))_{j \in \mathbb{N}}$ injective.

\begin{lemma}
Suppose that $\phi$ is a computable state on $\boldA$, and let $H^\#$ be a presentation of $H_\phi$ induced by the c.e. presentation $\boldA^\#$.  Then $H^\#$ and $\eta_\phi$ are computable.  In particular, then the cyclic vector $\eta_\phi(1_\boldA)$ is a computable point of $H^\#$, uniformly in a unital index for $\boldA^\#$.
\end{lemma}
\begin{proof}
The computability of $\eta_\phi$ is immediate from the definition.  To see that $H^\#$ is computable, given $j, k \in \mathbb{N}$, we have $\langle \eta_\phi(a_j), \eta_\phi(a_k) \rangle_\phi = \phi(a_k^*a_j)$, which is computable (uniformly in $j$ and $k$) because $\phi$ is computable.
\end{proof}

Given $a \in \boldA$, let $\pi_\phi(a) : \boldA/N_\phi \to \boldA/N_\phi$ be the linear operator defined by $\pi_\phi(a)(b+N_\phi) = ab+N_\phi$.  It is known that $\pi_\phi$ is contractive, so $\pi_\phi$ extends to a linear operator $\pi_\phi : H_\phi \to H_\phi$.  The GNS representation is the *-homomorphism $A \to B(H_\phi)$ given by $a \mapsto \pi_\phi(a)$.  When $H_\phi$ is infinite-dimensional, $B(H_\phi)$ is a non-separable \cstar-algebra, and hence cannot be given a computable presentation in the usual sense.  We therefore cannot ask that the GNS representation itself be a computable map.  Instead, the appropriate notion of computability is that $\pi_\phi(a)$ is a computable map of $H^\#$, uniformly in $a$.

\begin{theorem}
Suppose that $\phi$ is a computable state.  Let $H^\#$ be a presentation of $H_\phi$ induced by the c.e. presentation $\boldA^\#$.  If $a$ is a computable point of $\boldA^\#$ then $\pi_\phi(a)$ is a computable operator of $H^\#$.  Moreover, this is uniform in an $\boldA^\#$-name for $a$.
\end{theorem}
\begin{proof}
Suppose that $(a_j)_{j \in \mathbb{N}}$ is the effective enumeration of the generated points of $\boldA^\#$ that was used in defining $H^\#$.

Since $\boldA^\#$ is c.e., from a name for $a$ we can compute a rational number $M$ such that $\norm{a} < M$.  Suppose that we are given $n \in \mathbb{N}$ and an $H^\#$-name of a computable vector $x \in H$.  We can then compute $j$ such that $\norm{x - \eta_\phi(a_j)}_\phi < \frac{1}{M2^{n+1}}$.  Again using that $\boldA^\#$ is c.e., we can compute a rational number $M'$ such that $\norm{a_j} < M'$, and then we can compute $k$ such that $\norm{a-a_k} < \frac{1}{M'2^{n+1}}$.  Multiplication is computable on the generated points of $\boldA^\#$, so we can compute an $\ell$ such that $a_\ell = a_ka_j$.  Then $\eta_\phi(a_\ell)$ is a special point of $H^\#$, so it suffices to show that $\norm{\eta_\phi(a_\ell) - \pi_\phi(a)x}_\phi < 2^{-n}$.  For that, we calculate:
\begin{align*}
\norm{\eta_\phi(a_\ell) - \pi_\phi(a)x}_\phi &= \norm{\eta_\phi(a_ka_j) - \pi_\phi(a)(x)}_\phi \\
&\leq \norm{\eta_\phi(a_ka_j) - \eta_\phi(aa_j)}_\phi + \norm{\eta_\phi(aa_j) - \pi_\phi(a)(x)}_\phi \\
&\leq \norm{(a_k-a)}\norm{a_j} + \norm{\pi_\phi(a)(a_j) - \pi_\phi(a)(x)}_\phi\\
&< \frac{1}{M'2^{n+1}}M' + \norm{\pi_\phi(a)(a_j-x)}_\phi \\
&\leq 2^{-n-1} + \norm{a}\norm{a_j-x}_\phi \\
&< 2^{-n-1} + M\frac{1}{M2^{n+1}} \\
&< 2^{-n}
\end{align*}
\end{proof}

\section{Tracial states}\label{sec:Trace}
The trace simplex $T(\boldA)$, along with the canonical pairing between $T(\boldA)$ and $K_0(\boldA)$, plays a crucial role in the Elliott classification program.  It is therefore natural to hope that when $\boldA^\#$ is computable then $T(\boldA)$ can be given a computably compact evaluative presentation.  As we will see, this is not always the case.  We begin with some positive observations.  It will be useful to have the following reformulation of when a state is tracial, which follows immediately from the continuity of states and the polarization identity.

\begin{fact}
Let $\phi$ be a state and let $D \subseteq \boldA$ be a dense subset of $\boldA$.  Then $\phi$ is tracial if and only if $\phi(a^*a) = \phi(aa^*)$ for all $a \in D$.
\end{fact}

\begin{proposition}\label{lem:TraceComputable}
The trace space $T(\boldA)$ is a co-c.e. closed set in any evaluative presentation of $S(\boldA)$.
\end{proposition}
\begin{proof}
Let the generated points of $\boldA^\#$ be $(a_n)_{n \in \mathbb{N}}$.
For each $n \in \mathbb{N}$, let $e_n : S(\boldA) \to \mathbb{C}$ be defined by $e_n(\phi) = \abs{\phi(a_na_n^*)-\phi(a_n^*a_n)}$.  By definition of $S(\boldA)^\#$ being evaluative, each map $e_n$ is computable, uniformly in $n$.  Thus the sets $C_n = \{\phi \in S(\boldA) : e_n(\phi) = 0\}$ are co-c.e. closed, uniformly in $n$.  We have $T(\boldA) = \bigcap_{n \in \mathbb{N}}C_n$, which is thus a co-c.e. closed set of $S(\boldA)^\#$.
\end{proof}

In a few cases we can easily prove that $T(\boldA)$ inherits a computably compact evaluative presentation from $S(\boldA^\#)$.

\begin{proposition}\label{prop:TraceEvaluCompact}
If $\boldA$ is commutative or $\boldA$ has a unique trace then $T(\boldA)$ has a computably compact presentation evaluative over $\boldA^\#$.
\end{proposition}
\begin{proof}
If $\boldA$ is commutative then every state on $\boldA$ is tracial, so the result follows by Theorem \ref{thm:PresentationOfStates}.

If $\boldA$ has a unique trace then that trace is (uniformly) a computable point of $S(\boldA^\#)$ by \cite[Fact 3.25]{ComputablyCompact}.  In particular, $T(\boldA^\#)$ has a uniformly computable dense sequence of points, so by \cite[Lemma 3.27]{ComputablyCompact} it is c.e. closed in $S(\boldA^\#)$.  By Proposition \ref{lem:TraceComputable} and \cite[Proposition 3.29]{ComputablyCompact}, $T(\boldA)$ inherits a computably compact evaluative presentation from $S(\boldA^\#)$.
\end{proof}

In the Elliott invariant, traces are applied not directly to elements of $\boldA$, but to classes in $K_0(\boldA)$.  It is therefore useful to be able to effectively translate between projections in matrix algebras over $\boldA$ and classes in $K_0(\boldA)$.

\begin{definition}
A presentation $K_0(\boldA)^\#$ of $K_0(\boldA)$ is \emph{compatible with $\boldA^\#$} if there is an algorithm that takes as input an index of a projection $p \in M_n(\boldA)$ and outputs a $K_0(\boldA)^\#$-label for $[p]_{K_0}$ and there is an algorithm that takes as input a $K_0(\boldA)^\#$-label for a class $x \in K_0(\boldA)$ and outputs indices for projections $p \in M_n(\boldA)$ and $q \in M_k(\boldA)$ such that $x = [p]_{K_0} - [q]_{K_0}$.
\end{definition}

In \cite{UHFPaper} it was shown that there is a computable functor taking c.e. presentations of \cstar-algebras to c.e. presentations of their (ordered) $K_0$ groups.  We denote the presentation of $K_0(\boldA)$ thus obtained from $\boldA^\#$ by $K_0(\boldA^\#)$.  The only thing we will need about the construction of $K_0(\boldA^\#)$ is Theorem \ref{fact:FindProjection} below.  A proof of Theorem \ref{fact:FindProjection} is implicit in \cite{UHFPaper}, but we take this opportunity to provide it in detail since the result itself is of independent interest.  For all of the undefined notation used in the proof, we refer the reader to \cite[Section 2]{UHFPaper}.

\begin{theorem}\label{fact:FindProjection}
Up to computable isomorphism, the only presentation of $K_0(\boldA)$ that is compatible with $\boldA^\#$ is $K_0(\boldA^\#)$.
\end{theorem}
\begin{proof}
We first show that $K_0(\boldA^\#)$ is compatible with $\boldA^\#$.
Suppose we are given a $K_0(\boldA^\#)$-label $w$ for an element $x \in K_0(\boldA)$.  The construction of $K_0(\boldA^\#)$ in \cite[Proposition 2.14]{UHFPaper} arranges that the canonical surjection $\pi : \mathcal{D}(\boldA^\#) \times \mathcal{D}(\boldA^\#) \to K_0(\boldA^\#)$ is computable and the Grothendieck equivalence relation $\sim$ on $(\mathcal{D}(\boldA) \times \mathcal{D}(\boldA))^\#$ is c.e., all uniformly in an index for $\boldA^\#$.  Thus, by enumerating $\sim$, we can search for a pair $(g_0, g_1)$ of $\mathcal{D}(\boldA^\#)$ labels for $[p']_{MvN}, [q']_{MvN} \in \mathcal{D}(\boldA)$ such that $x = [p']_{K_0} - [q']_{K_0}$.  Since $\mathcal{D}(\boldA^\#)$ is computably supported by $\boldA^\#$, uniformly in $\boldA^\#$ (see \cite[Theorem 2.9]{UHFPaper}), from $g_0$ it is possible to compute $n$ and an index for a computable projection $p \in M_n(\boldA)$ such that $p \sim_{MvN} p'$, and similarly from $g_1$ we can compute $k$ and an index for a computable projection $q \in M_k(\boldA)$ such that $q \sim_{MvN} q'$.  Then $x = [p']_{K_0}-[q']_{K_0} = [p]_{K_0}-[q]_{K_0}$.

Now suppose that we have an index for a computable projection $p \in M_n(\boldA)$.  Let $(p_{m,k})_{m, k \in \mathbb{N}}$ be the sequence of projections used in the construction of $\mathcal{D}(\boldA^\#)$ in the proof of \cite[Theorem 2.9]{UHFPaper}, so $(p_{m,k})$ is a computable sequence of projections in $M_m(\boldA^\#)$, uniformly in $m$, and each projection in $M_m(\boldA)$ is Murray-von Neumann equivalent to some $p_{m,k}$.  By \cite[Corollary 2.2]{UHFPaper}, there is a c.e. subset $Q \subseteq \mathbb{N}\times\mathbb{N}$ such that if $e$ and $e'$ are indices of projections then $(e, e') \in Q$ if and only if these two projections are Murray-von Neumann equivalent, and an index of $Q$ can be computed from an index of $\boldA^\#$.  Since $p$ is Murray-von Neumann equivalent to some $p_{n,k}$, by searching $Q$ we can compute a value of $k$ for which this is the case.  Then, by the definition of $\mathcal{D}(\boldA^\#)$, $\langle n, k\rangle$ is a $\mathcal{D}(\boldA^\#)$-label for $[p]_{MvN}$.  Finally, by \cite[Proposition 2.14]{UHFPaper}, we can compute an index for $K_0(\boldA^\#)$ and an index for the canonical map $\gamma : \mathcal{D}(\boldA) \to K_0(\boldA)$.  Thus, by definition of $\gamma$ being a computable map (\cite[Definition 2.5]{UHFPaper}), from the $\mathcal{D}(\boldA^\#)$-label $\langle n, k\rangle$ for $[p]_{MvN}$ we can compute a $K_0(\boldA)$ label for $\gamma([p]_{MvN})$, i.e., for $[p]_{K_0}$.

We now turn to showing uniqueness.  Suppose that $K_0(\boldA)^\dagger$ is another presentation of $K_0(\boldA)$ that is compatible with $\boldA^\#$.  Given a $K_0(\boldA)^\dagger$-label for $x \in K_0(\boldA)$, using that $K_0(\boldA)^\dagger$ is compatible with $\boldA^\#$ we can compute $\boldA^\#$-indices for projections $p \in M_n(\boldA)$ and $q \in M_k(\boldA)$ such that $x = [p]_{K_0}- [q]_{K_0}$.  Then, using that $K_0(\boldA^\#)$ is compatible with $\boldA^\#$, from the indices for $p$ and $q$ we can compute $K_0(\boldA^\#)$-labels for $[p]_{K_0}$ and $[q]_{K_0}$, and hence we can compute a $K_0(\boldA^\#)$-label for $x = [p]_{K_0}-[q]_{K_0}$.  Thus the identity map is computable from $K_0(\boldA)^\dagger$ to $K_0(\boldA^\#)$.  Similarly, the identity map is computable from $K_0(\boldA^\#)$ to $K_0(\boldA)^\dagger$, so $K_0(\boldA)^\dagger$ is computably isomorphic to $K_0(\boldA^\#)$.
\end{proof}

A trace $\tau \in T(\boldA)$ extends naturally to a trace $\tau_n$ on $M_n(\boldA)$ by 
\[\tau_n([a_{i,j}]_{1 \leq i, j \leq n}) = \sum_{i=1}^{n}\tau(a_{i,i}).\]  Given a class $x \in K_0(\boldA)$, if $x = [p]_{K_0} - [q]_{K_0}$ for some projections $p \in M_n(\boldA)$ and $q \in M_k(\boldA)$, the evaluation of $\tau$ at $x$ is defined by $\tau(x) = \tau_n(p) - \tau_k(q)$.  It is well-known that this definition is independent of the projections $p$ and $q$ used to represent $x$.  The pairing of $T(\boldA)$ with $K_0(\boldA)$ is then defined by
\[\langle \tau, x \rangle = \tau(x).\]
\begin{proposition}\label{prop:EvaluateTrace}
If $T(\boldA)^\#$ is evaluative over $\boldA^\#$ then the pairing map is a computable map from the presentation of $T(\boldA) \times K_0(\boldA)$ induced by $T(\boldA)^\#$ and $K_0(\boldA^\#)$.
\end{proposition}
\begin{proof}
Suppose we are given a $K_0(\boldA^\#)$-label for $x \in K_0(\boldA)$, and a computable point $\tau \in T(\boldA^\#)$.  Theorem \ref{fact:FindProjection} lets us compute indices for projections $p \in M_n(\boldA)$ and $q \in M_k(\boldA)$ such that $x = [p]_{K_0}-[q]_{K_0}$.  Since $T(\boldA)^\#$ is evaluative, we can then compute $\tau_n(p)$ and $\tau_k(q)$.
\end{proof}

We note that the proof above is uniform, in the sense that from indices for $\boldA^\#$ and the evaluation map $T(\boldA) \times \boldA \to \mathbb{C}$ it computes an index for the trace pairing map.

Given $x \in K_0(\boldA)$, let $I_x = \{\tau(x) : \tau \in T(\boldA)\} \subseteq \mathbb{R}$ be the range of the trace pairing on $x$.  Since $T(\boldA)$ is compact and convex, and the pairing is continuous in the first coordinate, $I_x$ is a compact interval (provided $T(\boldA) \neq \emptyset)$.

\begin{corollary}\label{cor:ComputeTraceImage}
Suppose that $T(\boldA) \neq \emptyset$ and $T(\boldA)^\#$ is computably compact and evaluative over $\boldA^\#$.  For every $x \in K_0(\boldA)$, the interval $I_x$ is computably compact.  Moreover, from a $K_0(\boldA^\#)$-label for $x$ it is possible to compute the endpoints of $I_x$.
\end{corollary}
\begin{proof}
It suffices to show the latter claim, as a compact interval with computable endpoints is computably compact.  Let $f : T(\boldA) \to \mathbb{R}$ be defined by $f(\tau) = \tau(x) = \langle \tau, x \rangle$.  By Proposition \ref{prop:EvaluateTrace}, $f$ is computable, and an index for $f$ can be computed from an index for $\boldA^\#$ and an index for the evaluation map.  From an index for $f$ and an index for $T(\boldA^\#)$, since $T(\boldA^\#)$ is computably compact, we can compute $\sup_{\tau \in T(\boldA)}f(\tau)$ and $\inf_{\tau \in T(\boldA)}f(\tau)$ by \cite[Proposition 3.5]{ComputablyCompact}.
\end{proof}

We are now prepared to show that Proposition \ref{prop:TraceEvaluCompact} does not extend to all computable presentations of \cstar-algebras, and in fact not even to all computable presentations of simple AF algebras.

\begin{theorem}\label{thm:TraceNonComputable}
There is a simple AF algebra $\boldA$ and a computable presentation $\boldA^\#$ of $\boldA$ such that there is no computably compact presentation of $T(\boldA)$ that is evaluative over $\boldA^\#$.  Moreover, in this example no extreme trace of $\boldA$ is a computable map of $\boldA^\#$.
\end{theorem}
\begin{proof}
We make use of the Bratteli diagram constructed in \cite[Theorem 6.1]{Coronel}.  We recall the construction here.

Let $\beta \in (0, 1)$ be a right-c.e. non-computable real number, and fix a computable decreasing sequence of rationals $(b_n)_{n \in \mathbb{N}} \searrow \beta$ with $b_0 = 1$.  Let $\beta_0 = b_0 = 1$ and for each $n > 0$, let $\beta_n = b_n/b_{n-1}$.  Then $(\beta_n)_{n \in \mathbb{N}}$ is a computable sequence of rationals with $\prod_{n \in \mathbb{N}} \beta_n = \beta$.  Write $\frac{1+\beta_n}{2} = \frac{p_n}{q_n}$ in lowest terms and set $X_n = p_n$ and $Y_n = q_n-p_n$.  The Bratteli diagram $B$ has a single root level $L_0 = \{r\}$ and then, for every $n > 0$, the levels have two elements, $L_n = \{c_n, d_n\}$.  There is one edge from $r$ to $c_1$ and one edge from $r$ to $d_1$.  For $n \geq 1$, there are $X_n$ edges from $c_n$ to $c_{n+1}$, $X_n$ edges from $d_n$ to $d_{n+1}$, $Y_n$ edges from $c_n$ to $d_{n+1}$, and $Y_n$ edges from $d_n$ to $c_{n+1}$.  Let $\boldA$ be the AF algebra associated to $B$.  

Observe that $B$ is a computable Bratteli diagram, so by \cite[Main Theorem 3]{AFPaper} we obtain a computable presentation $\boldA^\#$ of $\boldA$.  Moreover, in this presentation, the projection $p \in \boldA$ corresponding to the vertex $c_1$ of the Bratteli diagram is a computable point.  By Corollary \ref{cor:ComputeTraceImage}, the maximum of $I_p = \{\tau(p) : \tau \in T(\boldA)\}$ is computable.  As noted in the proof of \cite[Theorem 6.1]{Coronel}, we have $I_p = \left[\frac{1-\beta}{2}, \frac{1+\beta}{2}\right]$.  Thus $\frac{1+\beta}{2}$ is a computable real number, contradicting the choice of $\beta$.

The final statement about no extreme trace being computable holds because in this algebra there are exactly two extreme traces and these traces are where the endpoints of $I_p$ are attained.
\end{proof}

Our final goal is to show that, for the presentation $\boldA^\#$ in Theorem \ref{thm:TraceNonComputable}, there is no computably compact presentation of $T(\boldA)$ such that the trace pairing with $K_0(\boldA^\#)$ is computable.  In general, knowing that $\langle \tau, \cdot \rangle = \langle \tau', \cdot \rangle$ does not imply that $\tau = \tau'$, so we do not expect that, in general, being able to compute $\langle \cdot, \cdot \rangle$ implies that the presentation of $T(\boldA)$ is evaluative.  However, if $\boldA$ is real rank $0$ then knowing the value of a trace on $K_0$ classes is sufficient to recover the trace itself.  Towards showing an effective version of this fact, we define the following effective version of real rank $0$.

\begin{definition}
$\boldA^\#$ is \emph{effectively real rank $0$} if there is an algorithm that takes as input a computable self-adjoint $a \in \boldA$ and an $n \in \mathbb{N}$ and outputs $k \in \mathbb{N}$, $\lambda_1, \ldots, \lambda_k \in \mathbb{Q}$, and indices for mutually orthogonal computable projections $p_1, \ldots, p_k \in \boldA$ such that $\norm{a - \sum_{j=1}^k\lambda_jp_j} < 2^{-n}$.
\end{definition}

Since the computable self-adjoint elements are dense in the set of all self-adjoint elements of $\boldA$, it follows that if $\boldA^\#$ is effectively real rank $0$ then $\boldA$ is real rank $0$.

\begin{proposition}\label{prop:RRZTrace}
If $T(\boldA)^\#$ is a presentation of $T(\boldA)$ such that the trace pairing is computable, and if $\boldA^\#$ is effectively real rank $0$, then $T(\boldA)^\#$ is evaluative over $\boldA^\#$.
\end{proposition}
\begin{proof}
We must show that $\operatorname{ev} : T(\boldA)^\# \times \boldA^\# \to \mathbb{C}$ is computable.  Suppose that we are given a name for $\tau \in T(\boldA)$ and a name for $x \in \boldA$. 

Let $a = \frac{1}{2}(x+x^*)$ and $b = \frac{1}{2i}(x-x^*)$.  Then $a=a^*$, $b=b^*$, $x = a+ib$, and $a$ and $b$ are computable from $x$.  Since $\tau$ is linear, this shows that it suffices to show that $\op{ev}$ is computable when the input element of $\boldA$ is self-adjoint.  We thus assume $x=x^*$.

Since $\boldA^\#$ is effectively real rank $0$, for each $n \in \mathbb{N}$ we can compute $k_n \in \mathbb{N}$, $\lambda_{n,1}, \ldots, \lambda_{n,k_n} \in \mathbb{Q}$, and indices for mutually orthogonal projections $p_{n,1}, \ldots, p_{n, k_n}$ such that $\norm{x - \sum_{j=1}^{k_n}\lambda_{n,j}p_{n,j}} < 2^{-n}$.  Let $a_n = \sum_{j=1}^{k_n}\lambda_{n,j}p_{n,j}$.  Using Theorem \ref{fact:FindProjection}, we can compute a $K_0(\boldA^\#)$-label for $[p_{n,j}]_{K_0}$ for each $n \in \mathbb{N}$ and each $j \in \{1, \ldots, k_n\}$, uniformly in $n$ and $j$.  We then have:
\begin{align*}
\operatorname{ev}(\tau, a_n) &= \tau(a_n) \\
&= \sum_{j=1}^{k_n}\lambda_{n,j}\tau(p_{n,j}) \\
&= \sum_{j=1}^{k_n}\lambda_{n,j}\langle\tau, [p_{n,j}]_{K_0}\rangle
\end{align*}
Since the trace pairing with $K_0$ is computable, this shows that $(\op{ev}(\tau, a_n))_{n \in \mathbb{N}}$ is a computable sequence.  Since $\tau$ is contractive, we have
\[\abs{\op{ev}(\tau, a_n) - \op{ev}(\tau, x)} = \abs{\tau(a_n) - \tau(x)} \leq \norm{a_n - x} < 2^{-n}.\]
Therefore $(\op{ev}(\tau, a_n))_{n \in \mathbb{N}}$ is a computable name for $\op{ev}(\tau, x)$.
\end{proof}

We do not expect that all computable presentations of real rank $0$ \cstar-algebras are effectively real rank $0$.  For our current goal of showing that the Elliott invariant is not always fully computable, the following class of examples is sufficient.

\begin{proposition}\label{prop:AFRRZ}
Every computable presentation of an AF algebra is effectively real rank $0$.
\end{proposition}
\begin{proof}
Let $a$ be a computable self-adjoint element of $\boldA$ and fix $n \in \mathbb{N}$.  By \cite[Lemma 7.2]{AFPaper} we can compute an index for a c.e. closed finite-dimensional *-subalgebra $F$ of $\boldA$ with $1_\boldA \in F$ and $d(a,F)<\frac{2^{-n}}{4}$; moreover, the proof of \cite[Lemma 7.2]{AFPaper} in fact computes a set of matrix units for $F$, so we can also compute $\dim(F)$.  Using \cite[Lemma 7.1]{AFPaper}, we can compute a finite multi-set $J$ and an isomorphism $\psi : \bigoplus_{j \in J}M_j(\mathbb{C}) \to F$.

Since $d(a,F) < \frac{2^{-n}}{4}$, we can search to find a $b \in F$ with $\norm{a - b} < \frac{2^{-n}}{4}$.  Observe that $\left(\frac{1}{2}(b+b^*)\right)^* = \frac{1}{2}(b+b^*)$, $\frac{1}{2}(b+b^*) \in F$, and 
\begin{align*}
\norm{a-\frac{1}{2}(b+b^*)}
&= \norm{\frac{1}{2}((a-b)+(a-b^*))} \\
&\leq \frac{1}{2}(\norm{a-b}+\norm{a-b^*}) \\
&= \frac{1}{2}(\norm{a-b}+\norm{(a-b)^*}) \\
&= \norm{a-b}
\end{align*}

Replacing $b$ with $\frac{1}{2}(b+b^*)$, we may therefore assume that $b=b^*$.

Let $c = \psi^{-1}(b)$, say $c = \oplus_{j \in J}c_j$ with each $c_j \in M_j(\mathbb{C})$.  For each $j$, search for $d_j \in M_j(\mathbb{Q}(i))$ such that $\norm{c_j-d_j} < \frac{2^{-n}}{4}$.  As above, replacing $d_j$ by $\frac{1}{2}(d_j+d_j^*)$ allows us to assume that $d_j$ is self-adjoint without increasing $\norm{c_j-d_j}$.  The self-adjoint matrices with distinct eigenvalues are dense in the set of self-adjoint elements of $M_j(\mathbb{Q}(i))$, and for a matrix with entries in $\mathbb{Q}(i)$ we can test whether the eigenvalues are distinct by computing the discriminant of the characteristic polynomial (which has coefficients in $\mathbb{Q}(i)$) and checking if it is $0$.  Thus, for each $j$, we can search for $e_j \in M_j(\mathbb{Q}(i))$ such that $e_j = e_j^*$, the eigenvalues of $e_j$ are distinct, and $\norm{d_j-e_j} < \frac{2^{-n}}{4}$.  In particular, we know the number of distinct eigenvalues of $e_j$ is exactly $j$, so by \cite[Theorem 2]{Brattka} we can effectively diagonalize $e_j$; that is, we can compute the distinct $\mu_{j, 1}, \ldots, \mu_{j,j} \in \mathbb{C}$ and indices for mutually orthogonal projections $q_{j,1}, \ldots, q_{j,j} \in M_j(\mathbb{C})$ such that $e_j = \sum_{k=1}^j \mu_{j,k}q_{j,k}$.  For each $j$ and $k$, compute $\lambda_{j,k} \in \mathbb{Q}$ such that $\abs{\lambda_{j,k} - \mu_{j,k}} < \frac{2^{-n}}{4}$.  For each $j$ let $f_j = \sum_{k=1}^j \lambda_{j,k}q_{j,k}$.  Note that, for each $j$,
\[\norm{e_j-f_j} = \norm{\sum_{k=1}^j(\mu_{j,k}-\lambda_{j,k})q_{j,k}} = \max_k \abs{\mu_{j,k} - \lambda_{j,k}} < \frac{2^{-n}}{4}.\]
Let $f = \oplus_{j\in J}f_j$.  After reindexing, we have computed $\lambda_1, \ldots, \lambda_r \in \mathbb{Q}$ and mutually orthogonal projections $q_1, \ldots, q_r \in \bigoplus_{j\in J}M_j(\mathbb{C})$ such that $f = \sum_{k=1}^r \lambda_kq_k$.

For each $k$ let $p_k = \psi(q_k)$.  Then the $p_k$ are mutually orthogonal computable projections in $F$, hence also in $\boldA$, and
\begin{align*}
\norm{a - \sum_{k=1}^r\lambda_kp_k} &\leq \norm{a - b} + \norm{b - \sum_{k=1}^r\lambda_kp_k} \\
&< \frac{2^{-n}}{4} + \norm{\psi^{-1}\left(b - \sum_{k=1}^r\lambda_kp_k\right)} \\
&= \frac{2^{-n}}{4} + \norm{c - f} \\
&\leq \frac{2^{-n}}{4} + \norm{c-d} + \norm{d-e} + \norm{e-f} \\
&= \frac{2^{-n}}{4} + \norm{\oplus_{j}(c_j-d_j)} + \norm{\oplus_{j}(d_j-e_j)} + \norm{\oplus_j(e_j-f_j)} \\
&= \frac{2^{-n}}{4} + \max_j\norm{c_j-d_j} + \max_j\norm{d_j-e_j}+\max_j\norm{e_j-f_j} \\
&< \frac{2^{-n}}{4} + \frac{2^{-n}}{4} + \frac{2^{-n}}{4} + \frac{2^{-n}}{4}\\
&= 2^{-n}
\end{align*}
\end{proof}

\begin{corollary}
There is a simple AF algebra $\boldA$ with a computable presentation $\boldA^\#$ such that $\op{Ell}(\boldA)$ has no fully computable presentation over $\boldA^\#$.
\end{corollary}
\begin{proof}
Let $\boldA$ and $\boldA^\#$ be as in the statement of Theorem \ref{thm:TraceNonComputable}.  Since $\boldA$ is AF, by Proposition \ref{prop:AFRRZ}, $\boldA^\#$ is effectively real rank $0$.  Thus by Proposition \ref{prop:RRZTrace}, if $\op{Ell}(\boldA)$ has a fully computable presentation over $\boldA^\#$, then $T(\boldA)$ has a computably compact presentation evaluative over $\boldA^\#$, contradicting Theorem \ref{thm:TraceNonComputable}.
\end{proof}

\section{Appendix (with Joaquin Cisneros): Computability of the unit}
In the context of unital \cstar-algebras, many calculations and descriptions of sets naturally involve expressions involving the unit $1_\boldA$ of $\boldA$, so in computability contexts it is often helpful to know that $1_\boldA$ is a computable point of the presentation.  For example, in Lemma \ref{lemma:StatesCoCE} above, to prove that the image of $S(\boldA)$ in a suitable product of discs is co-c.e. closed, we needed a computable name for $1_\boldA$.

In \cite[Theorem 4.1]{Gelfand} it was shown that the unit is always a computable point of every computable presentation of a commutative \cstar-algebra.  This was extended in \cite[Proposition 4.5]{UHFPaper} to the case where $\boldA$ is stably finite.  Subsequent work, such as earlier versions of \cite{FuncCalc} and \cite{EffectiveNuclearityPaper}, included the computability of the unit as a hypothesis when considering presentations of general unital \cstar-algebras.  Here, we show that the additional hypothesis is unnecessary.  In fact, this applies not only to computable presentations, but the broader class of (right-)c.e. presentations (that is, presentations where the norm is only required to be uniformly upper semicomputable, rather than uniformly computable, on generated points; see \cite[Definition 1.2]{UHFPaper}).

\begin{proposition}
The unit is a computable point of every c.e. presentation of a unital \cstar-algebra.
\end{proposition}
\begin{proof}
Since the generated points of $\boldA^\#$ are dense, there is a generated point $a$ such that $\norm{a-1_\boldA} \leq 1/2$.  Fix such an $a$.  Since $\norm{a-1_\boldA} < 1$, we have
\[a^{-1} = \sum_{n=0}^{\infty}(1_\boldA-a)^n,\]
so
\[\norm{a^{-1}} \leq \sum_{n=0}^{\infty}\norm{(1_\boldA-a)^n} \leq \sum_{n=0}^{\infty}2^{-n} = 2.\]
Note also that $\norm{a} \leq \norm{a-1_\boldA} + \norm{1_\boldA} \leq 2$.

Fix $n \in \mathbb{N}$.  Since $\boldA^\#$ is c.e., we can use an effective enumeration of the generated points and dovetailing to search for a generated point $c_n$ such that $\norm{c_na-a}<2^{-n-1}$.  To see that the search succeeds, note that since the generated points are dense there is a generated point $c$ with $\norm{c-1_\boldA} < 2^{-n-2}$, and for such a $c$, we have
\[\norm{ca-a} \leq \norm{c-1_\boldA}\norm{a} < 2\norm{c-1_\boldA} \leq 2^{-n-1}.\]
Then, having found a generated point $c_n$ such that $\norm{c_na-a} < 2^{-n-1}$, we have
\[\norm{c_n-1_\boldA} = \norm{(c_na-a)a^{-1}} \leq \norm{c_na-a}\norm{a^{-1}} \leq 2\norm{c_na-a} \leq 2^{-n}.\]
Thus $(c_n)_{n \in \mathbb{N}}$ is a computable name for $1_\boldA$.
\end{proof}

The proof above is non-uniform because of the need to fix $a$.  McNicholl \cite[Theorem 5.3]{EvaluativePresentations} has shown that a uniform proof is not possible; specifically, there is no uniform procedure that takes an index for a computable presentation $\boldA^\#$ of a unital \cstar-algebra and produces an $\boldA^\#$-index for $1_\boldA$.  

The concept of a \emph{unital index} for a c.e. presentation $\boldA^\#$ was introduced in \cite{FuncCalc}.  A unital index is $\langle e_1, e_2 \rangle$, where $e_1$ is an index for $\boldA^\#$, $e_2$ is an $\boldA^\#$-index for $1_\boldA$, and $\langle \cdot, \cdot \rangle$ is the standard computable bijection from $\mathbb{N}^2$ to $\mathbb{N}$.  The above observations can be summarized by saying that every c.e. presentation of a unital \cstar-algebra has a unital index, but a unital index cannot always be computed from just an index for the presentation.
\bibliographystyle{amsplain}
\bibliography{bib}
\end{document}